\documentclass[12pt,twoside]{amsart}
\usepackage{amssymb}

\nonstopmode

\numberwithin{equation}{section}
\font\tengothic=eufm10 scaled\magstep 1
\font\sevengothic=eufm7 scaled\magstep 1
\newfam\gothicfam
          \textfont\gothicfam=\tengothic
          \scriptfont\gothicfam=\sevengothic

\newcommand{\Z}{\mathbb{Z}}

\newcommand {\PP}{\mathbb{P}}
\newcommand {\CC}{\mathbb{C}}

\newcommand{\cI}{\mathcal{I}}

\newcommand{\cE}{\mathcal{E}}
\newcommand{\cO}{\mathcal{O}}
\DeclareMathOperator{\codim}{codim}
\DeclareMathOperator{\depth}{depth}

\DeclareMathOperator{\id}{id}

\DeclareMathOperator{\rank}{rank}

\newcommand{\ffi}{\varphi}
\newcommand{\al}{\alpha}

\DeclareMathOperator{\pnt}{\raise 0.5mm \hbox{\large\bf.}}

\DeclareMathOperator{\Proj}{Proj}

\newcommand{\lra}{\longrightarrow}
\newcommand{\im}{\operatorname{im}}

\newcommand{\mif}{\mbox{if} ~}

\newcommand{\HH}{H_{\mathfrak m}}

\newtheorem{theorem}{Theorem}[section]
\newtheorem{lemma}[theorem]{Lemma}
\newtheorem{proposition}[theorem]{Proposition}

\newtheorem{conjecture}[theorem]{Conjecture}

\newtheorem{problem}[theorem]{Problem}

\theoremstyle{definition}
\newtheorem{remark}[theorem]{Remark}
\newtheorem{example}[theorem]{Example}

\newtheorem{question}[theorem]{Question}

\begin{document}

\title[Resolutions of  subschemes of small degree]
{Minimal free resolutions of projective subschemes of small degree}

\author[Uwe Nagel]{Uwe Nagel}

\address{Department of Mathematics,
University of Kentucky, 715 Patterson Office Tower,
Lexington, KY 40506-0027, USA}

\email{uwenagel@ms.uky.edu}


\maketitle

\section{Introduction}

For an arbitrary projective subscheme $X \subset \PP^n$, it is
impossible to say much about its minimal free resolution. Some a
priori information is needed. In this article we will consider
integral subschemes as well as equidimensional Cohen-Macaulay
subschemes of small degree. There are several encouraging results
for such subschemes in the literature. Our goal is to discuss the
known results and to motivate further research in order to complete
these.

We first consider varieties, where we mean by a variety an
irreducible and reduced projective subscheme $X \subset \PP^n$.  We
will always consider $X$ as embedded in its linear span, that is, we
assume that $X$ is not contained in a hyperplane. In this case, we
say that $X \subset \PP^n$ is {\it non-degenerate}. Let $H \subset
\PP^n$ be a general hyperplane. Then Bertini's theorem implies that
$X \cap H \subset H \cong \PP^{n-1}$ is a reduced and non-degenerate
subscheme that is also irreducible provided $k := \dim X \geq 2$. It
follows that the intersection of $X$ with a general subspace of
dimension $n-k$ is a set of $d := \deg X$ points spanning a
$(n-k)$-dimensional linear space. This implies $d \geq n-k+1$, i.e.
$$
\deg X \geq \codim X + 1.
$$
Varieties with $\deg X = \codim X + 1$ are called {\it varieties of
  minimal degree}. They are completely classified. In particular, they
are always arithmetically Cohen-Macaulay and their minimal free
resolution is known. We recall some of the results in Section
\ref{sec-int}.

The next case, varieties satisfying $\deg X = \codim X + 2$ is
considerably more difficult. From a geometric viewpoint, Fujita has
a satisfying classification theory of these varieties (cf.\
\cite{F1}, \cite{F2}). Their cohomological properties are described
in \cite{HSV}. In particular, such varieties are not necessarily
arithmetically Cohen-Macaulay. The homological properties of them
are not well understood. More specifically, the authors posed in
\cite{HSV} the following problem.

\begin{problem} \label{prob-int}
Describe the minimal free graded resolution of a non-degenerate
variety $X \subset \PP^n$ with $\deg X = \codim X + 2$.
\end{problem}

Hoa was the first who attempted to solve this problem. In \cite{H}
he obtained partial information on the graded Betti numbers of the
varieties at hand. His results along with several  improvements are
presented in Section \ref{sec-int}. Our improvements are based on
the fact that every such curve $X$ is a divisor on a rational
surface scroll. This allows us to apply results in \cite{aBM-div}.
The best information, we get for varieties of codimension two (cf.\
Theorem \ref{thm-codim-2} below) and curves. However, in general, it
remains open to determine the graded Betti numbers of $X$. In some
sense, the first open case are curves in $\PP^5$ of degree six.
Remark \ref{rem-more-geom} suggests that a closer look at the
geometry of the curves is needed to make further progress.
\smallskip

Second, we drop the requirement on the subscheme $X \subset \PP^n$
to be integral. Then there is only  the trivial lower bound one for
its degree. Since the degree does not take lower dimensional
components into account, we assume that $X$ is equidimensional. Then
$X$ has degree one if and only if it is a linear space. Hence, the
first interesting case occurs when $\deg X = 2$. If $X$ is reducible 
then $X$ is simply the union of two linear spaces. However, if $X$
is irreducible, then its scheme structure can be rather complicated.
Nevertheless, for curves of degree two a complete description has
been obtained in \cite{NNS1} that allowed us to describe their
minimal free resolutions. This is described in Section
\ref{sec-equi}. In the case of space curves, there is also a  complete
classification of  curves having degree three \cite{Nollet} or four
\cite{Nollet-Schl}. A motivation for these results has been to
understand the Hilbert scheme of such curves, in particular whether
it is connected. In \cite{NNS1}, we also showed that double lines
can be used to construct curves with large cohomology.

It would be very interesting  to extend these results for curves to
subschemes of higher dimension and/or codimension.  Manolache has
obtained satisfying results for subschemes that have additional
properties like being a locally complete intersection. However, a
complete classification of Cohen-Macaulay subschemes of degree two
and three seems rather difficult. In fact, Vatne has shown that such
a classification should imply Hartshorne's conjecture on smooth
varieties of codimension two (cf.\ Remark \ref{rem-relation}).
\smallskip

While working on this paper it became clear that some results in the
literature can be improved. We  decided to include the stronger
statements and their proofs. This also illustrates  some of the
methods that are used to achieve these results.

It turns out that Problem \ref{prob-int}, as far as the Betti
numbers are concerned,  can be completely  solved for varieties of
codimension two. In fact, defining (by slight abuse of notation)
the depth of $X$ as the depth of its homogeneous coordinate ring we
show.

\begin{theorem} \label{thm-codim-2}
Let $X \subset \PP^n = \Proj (R)$ be a non-degenerate variety of codimension two and degree four. Then we have: 
\begin{itemize} 
\item[(a)] If $X$ is Cohen-Macaulay, then it is one of the varieties of the following list:
\begin{itemize}
\item[(i)] $X$ is a complete intersection cut out by two quadrics.
\item[(ii)] $X$ is a curve with minimal free resolution of the form
$$
0 \to R(-5) \to R(-4)^4 \to R(-3)^3 \oplus R(-2) \to I_X \to 0.
$$
\item[(iii)] $X$ is a surface of depth one and its minimal free
  resolution reads as 
$$
0 \to R(-6) \to R(-5)^5 \to R(-4)^{10} \to R(-3)^7  \to I_X \to 0.
$$
\end{itemize}
\item[(b)] In any case, the shape of the  minimal free resolution of $X$ is one  of the three types listed in {\rm (a)}. 
\end{itemize}
\end{theorem}

In particular we see that the dimension of such a variety is at most
two if it is Cohen-Macaulay, but not a complete intersection. 

While the varieties in case (ii) are divisors on a scroll, the surface
in case (iii) is not contained in any scroll. It is a generic projection of
the Veronese surface in $\PP^5$. 
\smallskip

Throughout this article, $X \subset \PP^n$ will be a projective
subscheme of dimension $k \geq 1$. We  denote its homogeneous
coordinate ring by $A = R/I_X$ where $I_X \subset R :=
K[x_0,\ldots,x_n]$ is the homogeneous ideal of $X$ and $K$ is a
field of arbitrary characteristic.

The $i$-th local cohomology module $\HH^i (A)$ of $A$ is isomorphic
to  $H^i_* (\cI_X) := \oplus_{j \in \Z} H^i (\cI_X (j))$, provided
$i \geq 1$.
\medskip

I am grateful to Irena Peeva for asking me to contribute this paper.
I enjoyed writing it. 

Moreover, I thank Nicolae Manolache, Scott Nollet, Roberto Notari, and Peter Schenzel for helpful comments. 

\section{Integral subschemes} \label{sec-int}

The focus of this section is Problem \ref{prob-int}. The solution of the analogous problem for varieties of minimal degree is known and we will make use of it.

Recall that a variety of minimal degree is either a quadric hypersurface, a rational normal scroll, or a cone over the Veronese surface in $\PP^5$. A modern proof of this classification can be found in \cite{EH}. Note that all varieties of minimal degree are arithmetically Cohen-Macaulay and that their minimal free resolutions are well-known. For simplicity, we just record their graded  Betti numbers here.

\begin{lemma} \label{lem-min-deg}
Let $X \subset \PP^n = \Proj (R)$ be a variety of minimal degree. Setting $c := \codim X$, $X$ has a minimal free resolution of the form
$$
0 \to R(-c-2)^{\alpha_c} \to \ldots \to R(-3)^{\alpha_2} \to R(-2)^{\alpha_1} \to I_X \to 0
$$
where ${\displaystyle \alpha_i = i \binom{c+1}{i+1}}$.
\end{lemma}

To see this, it suffices to consider rational normal curves. Their minimal free resolution is given by an Eagon-Northcott complex.

Among the varieties of minimal degree, the scrolls will be the most relevant here because Proposition 2.9 in \cite{aBM-div} says that effective integral divisors on a quadric of rank $ \geq 5$ and on cones over the Veronese surface in $\PP^5$ are always arithmetically Cohen-Macaulay. Following \cite{Ha} such scrolls can be described as follows.

Let $\cE = \oplus_{i=1}^k \cO_{\PP^1} (a_i)$ be a vector bundle over $\PP^1$ where $0 \leq a_1 \leq \ldots \leq a_k$ and $a_k > 0$. The tautological line bundle $\cO_{\PP(\cE)} (1)$ of its projectivization $\PP (\cE) \to \PP^1$ provides a birational map
$$
\PP (\cE) \to \PP^{k - 1 + \sum a_i}.
$$
We denote its image by $S(a_1,\ldots,a_k)$. It is a variety of dimension $k$ and minimal degree $c = \sum a_i$. Its homogeneous ideal is generated by the maximal minors of a $c \times 2$ matrix whose entries are linear forms.

Now we turn our attention to varieties $X$ satisfying
$$
\deg X = \codim X + 2.
$$
The cohomology of such varieties has been computed by Hoa, St\"uckrad, and Vogel.

\begin{theorem} \label{thm-coho}
If $X \subset \PP^n$ is a non-degenerate variety of degree $n-k+2 = \codim X + 2$. Then either $X$ is arithmetically Cohen-Macaulay, i.e.\ $\depth X = k+1$, or $1 \leq t := \depth X \leq k$ and
$$
H^i_*(\cI_X) \cong \left \{
\begin{array}{cl}
0 & \mif 1 \leq i \leq k \; \mbox{and} \; i \neq t \\
K[y_0,\ldots,y_{t-2}]^{\vee}(t-2)  & \mif i = t.
\end{array}
\right.
$$
where $y_0,\ldots,y_{t-2} \in A$ are algebraically independent elements of degree one and $N^{\vee}$ is the graded $K$-dual of $N$.
\end{theorem}

\begin{proof}
This follows from Theorems A and B in \cite{HSV}.
\end{proof}

\begin{remark}
Theorem \ref{thm-coho} implies that $X$ is arithmetically Buchsbaum, but not arithmetically Cohen-Macaulay if and only if $\depth X = 1$. Moreover, in this case we have $H^1_* (\cI_X) \cong K (-1)$.

Furthermore, if $t = \depth X \leq k$ then $X$ has the cohomology of a cone over a variety $Y \subset \PP^{n-t+1}$ of degree $n-k-t+3$ and depth one.
\end{remark}

In order to get information on the graded Betti numbers of these varieties we will induct on the dimension using general hyperplane sections.

\begin{remark}[First reduction] \label{rem-first-red}
Let $H \subset \PP^n$ be a general hyperplane defined by the linear form $l \in R$. Set $\overline{R} := R/l R$. Since the homogeneous coordinate ring $A$ of $X$ has positive depth, the graded Betti numbers of $A$ as $R$-module and the graded Betti numbers of $A/l A$ as $\overline{R}$-module coincide (cf., e.g., \cite{BH}, Proposition 1.1.4).

Moreover, if $\depth X = \depth A \geq 2$ then $A/l A$ is the homogeneous coordinate ring of $X \cap H \subset H \cong \PP^{n-1}$.
\end{remark}

Repeated use of this remark provides the graded Betti numbers if $X$ is  arithmetically Cohen-Macaulay.

\begin{proposition} \label{prop-acm}
Let $X \subset \PP^n$ be a non-degenerate integral variety of codimension $c$ and degree $c+2$ that is arithmetically Cohen-Macaulay. Then $X$ has a minimal free resolution of the form
\begin{eqnarray*}
\lefteqn{
0 \to R(-c-2) \to R(-c)^{\alpha_{c-1}} \to R(-c+1)^{\alpha_{c-2}}  \to } \\[1ex]
&& \hspace*{5cm} \ldots \to R(-2)^{\alpha_1} \to I_X \to 0
\end{eqnarray*}
where
\begin{eqnarray*}
\alpha_i  & = & i \binom{c+1}{i+1} - \binom{c}{i-1} \quad \mif \; 1 \leq i \leq c-1.
\end{eqnarray*}
\end{proposition}

\begin{proof}
This is Theorem 1 in \cite{H}. Using Lemma \ref{lem-curve-on-scroll} below, it could also be obtained from \cite{aBM-div}, Theorem 2.4.
\end{proof}

\begin{remark}
As observed in \cite{HSV}, Theorem B, the resolution above shows that $X$ is arithmetically Cohen-Macaulay if and only if it is arithmetically Gorenstein. This fact can be also deduced from a result of Stanley. Indeed, if $X$ is arithmetically Cohen-Macaulay then its h-vector must be $(1, c, 1)$. In particular, it is symmetric. Now \cite{S}, Theorem 4.4, shows that $X$ must be arithmetically Gorenstein.
\end{remark}

If we compute the graded Betti number of the variety $X$ that is not arithmetically Cohen-Macaulay, then using the first reduction (Remark \ref{rem-first-red}), we may assume that $\depth X = 1$. In this situation we will use.

\begin{remark}[Second reduction] \label{rem-second-red}
If $X$ has depth one, then we know by Theorem \ref{thm-coho} that $\HH^1 (A) \cong K(-1)$.
Hence, the long exact cohomology sequence associated to
$$
0 \to A(-1) \stackrel{l}{\longrightarrow} A \to A/l A \to 0
$$
implies $\HH^0 (A/ l A) \cong K(-2)$. This provides the exact sequence
$$
0 \to K(-2) \to A/ lA \to \overline{R}/I_{X \cap H} \to 0.
$$
Assuming by induction that we know the minimal free resolution of $X \cap H$,
the Horseshoe lemma gives a free resolution of $A/l A$ as $\overline{R}$-module which might not be minimal. But in any case,  we get upper bounds for the graded Betti numbers of $X$ by knowing upper bounds for the graded Betti numbers of $X \cap H$.
\end{remark}

In order to use this idea, we have to get information on curves first.  For curves, our starting point is a slight modification of  Proposition 1.2 in \cite{EEK} in order to include possibly singular curves.

\begin{lemma} \label{lem-curve-on-scroll}
If $C \subset \PP^n$ is a non-degenerate integral curve of degree $n+1$ that is not arithmetically Cohen-Macaulay, then $C$ is contained in a rational normal surface scroll.
\end{lemma}

\begin{proof}
Consider the exact sequence
$$
0 \to R/I_C \to H^0_*(\cO_C) \to H^1_*(\cI_C) \to 0.
$$
Since $H^1_*(\cI_C) \cong K(-1)$ by Theorem \ref{thm-coho}, this sequence becomes
$$
0 \to R/I_C \to H^0_*(\cO_C) \to K(-1) \to 0.
$$
It shows that $H^0_*(\cO_C)$ is a graded integral domain that is generated in degree one. Thus, it is the homogenous coordinate ring of an integral curve $C' \subset \PP^{n+1}$ of degree $n+1$, i.e.\ $C'$ is a rational normal curve. Moreover, $C$ is the projection of $C'$ from a point $P$ off the curve. Let $L$ be any line joining $P$ to a point of $C$. Then, by \cite{EEK}, Lemma 1.3, $C' \cup L$ lies on a  rational normal surface scroll, whose projection will be a scroll containing $C$.
\end{proof}

\begin{remark}
The analogous result is not true for integral curves of degree $n+2$ that span $\PP^n$ (cf.\ \cite{BS}).
\end{remark}

The lemma above says that each such curve is an effective divisor on a
scroll $S$. Thus, we can apply the results in \cite{aBM-div} where it is shown that the scroll must be smooth, and also the divisor class of $C$ is determined.
This provides.

\begin{proposition} \label{prop-betti-curve}
Let $C \subset \PP^n = \Proj (R)$ be a curve as in Lemma
\ref{lem-curve-on-scroll}. Then $I_C$ has a minimal free resolution of
the form
$$
0 \to F_n \to \ldots \to F_1 \to I_C \to 0
$$
where
$$
F_i = R(-i-1)^{\alpha_i + \beta_i} \oplus R(-i-2)^{\gamma_i}
$$
and
\begin{eqnarray*}
\alpha_i  & = & i \binom{n-1}{i+1}, \\[1ex]
\beta_i & = & \left \{ \begin{array}{cl}
n-3 & \mif \; i = 1 \\
0 & \mif \; i \geq n-2,
\end{array} \right.
\end{eqnarray*}
$$
\gamma_{n-3} = \binom{n+1}{3} - (n-1)^2,\;  \gamma_{n-2} = \binom{n}{2},\; \gamma_{n-1} = n+1,\;  \gamma_n = 1, \\[1ex]
$$
$$
\gamma_i - \beta_{i+1} = \binom{n+1}{i+1} - \binom{n-1}{i+1} - (n-1) \binom{n-2}{i} \quad \mif \; i \geq 1,
$$
and
$$
0 \leq \beta_i , \quad 0 \leq \gamma_i \leq \binom{n}{i} \quad \mif \; i \geq 1
$$
\end{proposition}

\begin{proof}
The estimate
$$
\gamma_i \leq \binom{n}{i}
$$
is given in \cite{H}, Proposition 8, the other claims follow from \cite{aBM-div}, Theorem 5.10. Hoa's inequality follows quickly from our second reduction principle. Indeed, by Remark \ref{rem-second-red} we have the exact sequence
$$
0 \to K(-2) \to A/ lA \to \overline{R}/I_{C \cap H} \to 0.
$$
Since the graded Betti numbers of the Gorenstein set $C \cap H$ of $n+1$ points are well known, an application of the Horseshoe lemma implies \cite{H}, Proposition 8.

Hoa's result can be refined using Lemma \ref{lem-curve-on-scroll}. In fact, this lemma  says that there is a rational normal surface scroll $S$ that contains $C$. Set $M := I_C/I_S$. Then the initial degree of $M$ is two and $\rank [M]_2 = n-3$. Hence, the length of the linear strand of the minimal free resolution of $M$ as $R$-module is at most $n-3$ by \cite{EK}, Theorem 1.1.  Now the Horseshoe lemma applied to
$$
0 \to I_S \to I_C \to M \to 0
$$
implies the remaining claims. For details we refer to the proof of Theorem 5.10 in \cite{aBM-div}.
\end{proof}

\begin{remark} \label{rem-compare-hoa}
Proposition \ref{prop-betti-curve} gives more precise information on
the Betti numbers of the curve  than \cite{H}, Proposition 8. The key
for the improvement is Lemma \ref{lem-curve-on-scroll}. 

Also observe that we do not use any assumption on the characteristic. 
\end{remark}

Let us illustrate our result for curves of low codimension. We have completely determined their graded Betti numbers if the codimension is at most three.

\begin{example} \label{ex-curves}
Let $C \subset \PP^n$ be a non-degenerate integral curve of degree $n+1$.

(i) If $C$ has codimension two then the shape of its minimal free resolution is either
$$
0 \to R(-4) \oplus R(-2)^2 \to I_C \to 0
$$
or
$$
0 \to R(-5) \to R(-4)^4 \to R(-3)^3 \oplus R(-2) \to I_C \to 0.
$$

(ii) If $C$ has codimension three then the shape of its minimal free resolution is either
$$
0 \to R(-5) \to   R(-3)^5 \to  R(-2)^5 \to I_C \to 0
$$
or
\begin{eqnarray*}
\lefteqn{
0 \to R(-6) \to R(-5)^5 \to R(-4)^6 \oplus R(-3)^2 \to } \\[1ex]
& & \hspace*{5cm} R(-3) \oplus R(-2)^4 \to I_C \to 0.
\end{eqnarray*}

(iii) The first case where we do not have complete information about the Betti numbers occurs when $C$ has codimension four. If $C$ is not arithmetically Cohen-Macaulay then Proposition \ref{prop-betti-curve} gives that $C \subset \PP^5$ has a minimal resolution of the form
\begin{eqnarray*}
\lefteqn{
0 \to R(-7) \to R(-6)^6 \to R(-5)^{10} \oplus R(-4)^3 \to } \\[1ex]
& & R(-4)^4 \oplus R(-3)^{8+\beta_2} \to R(-3)^{\gamma_1} \oplus R(-2)^8 \to I_C \to 0
\end{eqnarray*}
where $0 \leq \gamma_1 \leq 5$ and $\gamma_1 - \beta_2 = -3$, thus $\beta_2 \geq 3$.
\end{example}

The next example shows that in case $C$ is not arithmetically Cohen-Macaulay, there is not just one possible minimal free resolution, in general.

\begin{example} \label{ex-deg-6}
We consider two curves of degree six in $\PP^5$ that are not arithmetically Cohen-Macaulay.

(i) Let $C_1 \subset \PP^5$ be the curve with parametrization $\{s^6 : s^5 t : s^3 t^3 : s^2 t^4 : s t^5 : t^ 6\}$. Its minimal free resolution reads as
\begin{eqnarray*}
\lefteqn{
0 \to R(-7) \to R(-6)^6 \to R(-5)^{10} \oplus R(-4)^3 \to } \\[1ex]
& & R(-4)^4 \oplus R(-3)^{12} \to R(-3) \oplus R(-2)^8 \to I_{C_1} \to 0.
\end{eqnarray*}

(ii) Let $C_1 \subset \PP^5$ be the curve with parametrization $\{s^6 : s^5 t : s^4 t^2 : s^2 t^4 : s t^5 : t^ 6\}$. Its minimal free resolution has the form
\begin{eqnarray*}
\lefteqn{
0 \to R(-7) \to R(-6)^6 \to R(-5)^{10} \oplus R(-4)^3 \to } \\[1ex]
& & R(-4)^4 \oplus R(-3)^{11} \to  R(-2)^8 \to I_{C_2} \to 0.
\end{eqnarray*}

Hence, the ideal of $C_1$ has a minimal cubic generator whereas the ideal of $C_2$ is generated by quadrics.
\end{example}

It would be interesting to decide.

\begin{question}
Are the two resolutions described in Example \ref{ex-deg-6} the only possible resolutions of integral degree six curves in $\PP^5$ that are not arithmetically Cohen-Macaulay?
\end{question}

\begin{remark}
We expect that the larger the codimension of the curve, the more different sets of graded Betti numbers do occur.
\end{remark}

In order to shed some light onto this phenomenon we observe.

\begin{remark} \label{rem-more-geom}
Consider again the curves $C_1, C_2 \subset \PP^5$ of Example
\ref{ex-deg-6}. It is not difficult to see that the curve $C_1$ is
contained in a scroll $S(1, 3)$ whereas the curve $C_2$ is contained
in a scroll $S(2, 2)$, i.e.\ the two curves are contained in scrolls
of different type.

In general, Theorem 5.10 in \cite{aBM-div} implies that the divisor
class $H + 2 F$ of each scroll $S(a_0, a_1)$ with $1 \leq a_0 \leq a_1 =
n - 1 - a_0$ contains an integral curve $C \subset \PP^n$ of degree
$n+1$. This leaves  $\lfloor \frac{n-1}{2} \rfloor$
possible choices for the scrolls that contain such a curve. Example
\ref{ex-deg-6} suggests that the geometry of the scrolls plays a
crucial role in determining the minimal free resolution of $C$.
\end{remark}

Now we are ready to prove Theorem \ref{thm-codim-2}.

\begin{proof}[Proof of Theorem \ref{thm-codim-2}]  
Because of Proposition \ref{prop-acm} it remains to consider
non-arithmetically Cohen-Macaulay varieties.  For these, Theorem \ref{thm-coho} shows that $X$ is Cohen-Macaulay if and only if $X$ has depth one. Thus, for the remainder of the proof of (a), we may
assume that the depth of $X$ is one.

If $X$ is a curve then the results follows by  Example \ref{ex-curves}.

Now let $X$ be a surface of depth one. Then its general hyperplane
section is a curve $C$ of depth one, whose graded Betti numbers we
already know. Hence, the second reduction (Remark
\ref{rem-second-red}) provides the following exact diagram
$$
\begin{array}{ccccccccc}
0 \lra & K(-2) & \lra & A/l A & \lra & \overline{R}/I_C & \lra 0. \\
& \uparrow & & & & \uparrow \\
& \overline{R}(-2) &&&& \overline{R} \\
& \uparrow & & & & \uparrow \\
& \overline{R}(-3)^4 & && & \overline{R}(-2) \oplus \overline{R}
(-3)^3 \\
& \uparrow & & & & \uparrow \\
& \overline{R}(-4)^6 &&&& \overline{R}(-4)^4 \\
& \uparrow & & & & \uparrow \\
& \overline{R}(-5)^4 &&&& \overline{R}(-5) \\
& \uparrow & & & & \uparrow \\
& \overline{R}(-6) &&&& 0 \\
& \uparrow \\
& 0
\end{array}
$$
Applying the Horseshoe lemma, we get a non-minimal free resolution
of $A/ l A$. Since $A/ l A$ is a cyclic $\overline{R}$-module, we
can cancel the direct summands $\overline{R}(-2)$ at the beginning
of this resolution. That provides a linear resolution which is
certainly minimal. 

Now assume that $X$ has depth one and that its dimension is at least three. We argue by induction on the dimension that this is impossible. Indeed, let $Y$ be the general hyperplane section of $X$. Then $Y$ has also depth one. We use again the exact sequence (cf.\ Remark  
\ref{rem-second-red}) 
$$
0 \lra  K(-2)  \lra  A/l A  \lra  \overline{R}/I_Y  \lra 0.
$$
It implies, as in the case of a surface above, that $I_Y$ has exactly one more quadratic generator than $I_X$. In particular, it must have at least one such generator. But we have just seen that $[I_Y]_2 = 0$ if $Y$ is a surface. This contradiction finishes the proof of (a) 

Claim (b) follows by combining the first reduction (cf.\ Remark \ref{rem-first-red}) and (a). 
\end{proof}

In general, the method of the above proof provides only a weaker
result. Since arithmetically Cohen-Macaulay  varieties are covered
by Proposition \ref{prop-acm}, we exclude them now.

\begin{theorem} \label{prop-betti-gen}
Let $X \subset \PP^n = \Proj (R)$ be a non-degenerate variety of
degree $c+2 := \codim X + 2$ such  that $\depth X \leq \dim X$.  Set
$p = n + 1 - \depth X$.

Then $I_X$ has a minimal free resolution of
the form
$$
0 \to F_p \to \ldots \to F_1 \to I_X \to 0
$$
where
$$
F_i = R(-i-1)^{\delta_i} \oplus R(-i-2)^{\gamma_i}
$$
and
\begin{eqnarray*}
\delta_1 & = & \binom{c+2}{2} - p - 2, \\
\delta_i & = & 0 \quad \mif \; i \geq c, \\
\gamma_i & = &  \left \{ \begin{array}{cl}
\binom{p+1}{i+1} & \mif \; i \geq c \\[1ex]
\binom{p+1}{c} - (c+1) & \mif \; i = c-1,
\end{array} \right. \\
\gamma_i & \leq & \left \{ \begin{array}{cl}
\binom{p+1}{c-1} - c^2 & \mif \; i = c-2 \\[1ex]
\binom{p+1}{i+1} & \mif \; 1 \leq i \leq c-3,
\end{array} \right.
\end{eqnarray*}
and
$$
\gamma_i - \delta_{i+1} = \binom{p+1}{i+1} + \binom{c}{i+2} - (c+1) \binom{c+1}{i+1} \quad \mif \; i \geq 1.
$$ 
\end{theorem}

\begin{proof}
This is for the most part Theorem 2 in \cite{H}. Using our reduction steps we sketch an alternative proof that also provides our improvement of \cite{H}, Theorem 2.

Again, we may assume that $\depth X = 1$. If $X$ is a curve then we have
already shown the stronger Proposition \ref{prop-betti-curve}.

Let $\dim X \geq 2$ and let $Y = X \cap H$ be a general hyperplane
section of $X$. Then the induction hypothesis applies to $Y$. Let us
use $\gamma_i (X), \delta_i (X), \gamma_i (Y), \delta_i (X)$ to
denote the  graded Betti numbers of $X$ and $Y$, respectively. Since
the resolution of $K$ is given by the Koszul complex, the second
reduction  (Remark \ref{rem-second-red}) shows that $X$ has a
minimal free resolution of the form
$$
0 \to F_p \to \ldots \to F_1 \to I_X \to 0
$$
where
$$
F_i = R(-i-1)^{\delta_i} \oplus R(-i-2)^{\gamma_i}
$$
and
$$
\gamma_i (X) = \gamma_i (Y) + \binom{p}{i}.
$$
As in the proof of Theorem \ref{thm-codim-2}, we see that
cancelation occurs at the beginning, i.e.\ $I_X$ has one quadric
generator less than $I_Y$. For reason of degree there cannot be any
cancelation at the end. Thus, an easy computation finishes the
proof.
\end{proof}

\begin{remark}
Let $X$ be as in the preceding theorem. Put $k := \dim X$ and assume
that $X$ is Cohen-Macaulay, i.e.\ it has depth one. Then $I_X$ has
exactly $\binom{c+1}{2} - k - 1$ minimal generators of degree two.
The ideal of a  variety of minimal degree and codimension $c-1$ is
generated by $\binom{c}{2}$ quadrics. Hence, $X$ cannot be a divisor
on a variety of minimal degree if $k \geq c$.

On the other hand, if $k \leq c$ then such a variety does exist on a
rational normal scroll by Theorem 5.12 in \cite{aBM-div}.
\end{remark}

If $X$ is a divisor on a variety of minimal degree, then we have
stronger information on its Betti numbers than in Theorem
\ref{prop-betti-gen}.

\begin{proposition}
Let $X \subset \PP^n = \Proj (R)$ be a non-degenerate variety of
degree $c+2$, dimension $k = n-c$, and depth one that is a divisor
on a variety of minimal degree.  Then $I_X$ has a minimal free
resolution of the form
$$
0 \to F_n \to \ldots \to F_1 \to I_X \to 0
$$
where
$$
F_i = R(-i-1)^{\alpha_i + \beta_i} \oplus R(-i-2)^{\gamma_i}
$$
and
\begin{eqnarray*}
\alpha_i  & = & i \binom{c}{i+1}, \\[1ex]
\beta_i & = & \left \{ \begin{array}{cl}
c-k-1 & \mif \; i = 1 \\
0 & \mif \; i \geq c-k,
\end{array} \right. \\[1ex]
\gamma_i & = & \left \{ \begin{array}{cl} \binom{n+1}{i+1} - c
\binom{c}{i+1} -  \binom{c-1}{i} & \mif \; c-k-1 \leq i < c
\\[1ex]
\binom{n+1}{i+1} & \mif \; c  \leq i \leq n,
\end{array} \right.
\end{eqnarray*}
$$
\gamma_i - \beta_{i+1} = \binom{n+1}{i+1} - \binom{c}{i+1} - c
\binom{c-1}{i} \quad \mif \; i \geq 1,
$$
and
$$
0 \leq \beta_i , \quad 0 \leq \gamma_i \leq \binom{n+1}{i} \quad
\mif \; i \geq 1
$$
\end{proposition}

\begin{proof}
This is a special case of Theorem 5.10 in \cite{aBM-div}.
\end{proof}

It is also interesting to investigate varieties of slightly higher
degree. For smooth varieties, geometric classifications have been obtained, for example, in \cite{O}, \cite{I}. A study of the minimal free resolutions has been begun by Brodmann and Schenzel in the
case of curves. For details, we refer to  [Brodmann-Schenzel].

\section{Equidimensional Cohen-Macaulay subschemes} \label{sec-equi}

Now we will consider more general subschemes. Throughout this
section $X \subset \PP^n$ denotes a  non-degenerate  equidimensional
Cohen-Macaulay subscheme, i.e.\ its homogeneous ideal $I_X$ is
unmixed and its homogeneous coordinate ring  $A = R/I_X$ is locally
Cohen-Macaulay. Note that a projective subscheme $Y$ has these
properties if and only if its cohomology has finite length, i.e.\
the modules $H^i_* (\cI_Y)$ have finite length for all $i$ with $1
\leq i \leq \dim Y$.

If $X$ is not a linear space then its degree is at least two. Here,
we will focus on the case $\deg X = 2$.  Note that every
equidimensional subscheme of dimension one, a curve, is
Cohen-Macaulay. This is no longer true in higher dimension. In fact,
if $C$ is a curve that is not arithmetically Cohen-Macaulay, then
any cone over it is not Cohen-Macaulay. It turns out that the
Cohen-Macaulay requirement puts rather strong conditions on $X$ as
the description of such schemes becomes considerably simpler when
the dimension is at least two and the codimension is sufficiently
small. The largest variety of such schemes occurs in dimension one.
Fortunately, there is a complete description of the homogeneous
ideals of degree two curves that also allows us to describe their
minimal free resolutions, including the maps. We will describe these
results first. Then we discuss Manolache's results in higher
dimension and conclude with Vatne's equivalence of his conjecture
about multiple structures of degree two and three and Hartshorne's
conjecture about smooth codimension two varieties.

Before we turn to irreducible curves, we discuss the simpler case of
reducible subschemes. We exclude  the trivial case of a
hypersurface.

\begin{proposition} \label{prop-red-2}
Let $X \subset \PP^n$  be a reducible, equidimensional, non-degenerate
subscheme of degree two and  dimension $k$ with $1 \leq k < n$. Then
$X$ is Cohen-Macaulay if and only if  $n \geq 3$ is odd and $X$ is the
disjoint union of two disjoint linear subspaces of dimension $
\frac{n-1}{2}$. 

Moreover, in this case its minimal free resolution is the tensor
product of the resolutions of its two linear subspaces. In particular,
the resolution has the form 
$$
0 \to R(-n-1)^{\delta_n} \to \ldots \to R(-3)^{\delta_2} \to
R(-2)^{\delta_1} \to I_X \to 0 
$$
where $\delta_i = \binom{n+1}{i+1} - 2 \binom{q}{i+1}$,  $1 \leq i
\leq n$,  and $q := \frac{n+1}{2}$. 
\end{proposition}

\begin{proof}
By assumption on the degree, $X$ is the union of two $k$-dimensio\-nal
linear spaces $L, M$. Thus, we have the exact sequence 
$$
0 \to R/I_X \to R/I_L \oplus R/I_M \to R/(I_L + I_M) \to 0.
$$
Set $p := \dim R/(I_L + I_M)$. Comparing the Hilbert functions in
degree one, we get $n+ 1 + p = 2 (k+1)$, thus in particular $p < k$
because $k <n$. 
Since $L$ and $M$ are arithmetically  Cohen-Macaulay, the sequence implies
$$
H^{i}_* (\cI_X) \cong \left \{ \begin{array}{cl}
0 & \mif 1 \leq i \leq k, \; i \neq p+1 \\
\HH^p (R/(I_L + I_M)) & \mif i = p+1.
\end{array} \right.
$$
Hence $X$ is not Cohen-Macaulay if $p > 0$ because then the module
$\HH^p (R/(I_L + I_M))$ has not finite length.  If $p = 0$, i.e.\ $L$
and $M$ do not meet, then $X$ is Cohen-Macaulay and $k =
\frac{n-1}{2}$. Moreover, 
in this case the tensor product of the minimal free resolutions of
$I_L, I_M$ gives the minimal free resolution of $I_X = I_L \cap I_M
\cong I_L \otimes I_M$. Using the formula 
$$
\sum_{i=0}^j \binom{q}{i+1} \binom{q}{j-i+1} = \binom{2q}{j+2} - 2
\binom{q}{j+2}. 
$$
we get the claim about the Betti numbers of $X$.
\end{proof}

It remains to consider irreducible subschemes of degree two. In the
case of a curve $C$, the support is a line $L$, thus $C$ is often
called a double line. Its ideal $I_C$ contains $I_L^2$, thus $C$ is a
rope. It will be useful to consider other ropes as well. A {\it rope
  $C$ (supported on the line $L$)} is a curve such that $I_L^2 \subset
I_C \subset I_L$. A $d$-rope 
is just a rope of degree $d$. We may and will assume that the
supporting line $L$ is defined by the ideal $I_L =
(x_0,\ldots,x_{n-2})$.  In order to stress the particular role played
by 
the line, we denote the coordinate ring of $\PP^n$ by $R =
K[x_0,\ldots,x_r, t, u]$ where 
$r := n-2$. Then the coordinate ring of the line $L$ is
$S := K[t, u]$.

The homogeneous ideal of a rope has been characterized in  Theorem 2.4
of \cite{NNS1}. Here, the ideal generated by the $k \times k$ minors
of a matrix $M$ is denoted by $I_k (M)$. 

\begin{theorem} \label{charact} Let $C \subset \PP^n$ be a curve of
degree at most $n-1$.
Then the following conditions are equivalent:
\begin{itemize}
\item[1.] $ C $ is an $(n-k)$--rope supported on the line $ L;$
\item[2.] $ I_C = ((I_L)^2, [x_0, \dots, x_r] B) $ where the matrix
$ B $ gives a graded homomorphism $ \varphi_B : \oplus_{j=1}^k
  S(-\beta_j-1) \to 
S^{r+1}(-1) $ such that $ \codim(I_k(B)) = 2;$
\item[3.] $ I_C = ((I_L)^2, F_1, \dots, F_k) $ where
$ V(F_1, \dots, F_k) \subset \PP^n$ is a scheme of codimension $
k$ which contains $L$ and is smooth at the points of $ L.$
\end{itemize}
\end{theorem}

\begin{remark} \label{rem-ass-matrices}
According to \cite{NNS1}, Remark 2.7 and Proposition 3.1, every $(n-k)$--rope supported
on the  line $L$ is related to matrices $A, B$ with entries in $S$
such that there is an  exact sequence
\begin{equation*} 
0 \to \oplus_{j=1}^k S(-\beta_j-1)
\stackrel{B}{\longrightarrow} S(-1)
\stackrel{A}{\longrightarrow} \oplus_{i=0}^{r-k} S(\al_i-1) 
\to H^1_* (\cI_C) \to 0
\end{equation*}
where the ideals of maximal minors $ I_k(B) $ and $
I_{r+1-k}(A)$ have codimension 2. Note that $A^t$ is just the
syzygy matrix of $B^t$.
\end{remark}

\begin{example} \label{ex-ropes-ideal}
(i) If $n=2$, then the matrix $B$ must be of the form $B =
\begin{bmatrix}
 F \\
G \end{bmatrix}$ where $\{F, G\} \subset S$ is a regular sequence.
Hence, we recover the result  that the ideal of every
double line in $\PP^3$ is of the form
$$
(x_0, x_1)^2 + (x_0 F + x_1 G).
$$
For further results on double lines in $\PP^3$ we refer to \cite{mig} and \cite{M-05}. 

(ii) Using in case $n = 5$ the matrix $B :=
\begin{bmatrix}
 t & 0 & 0  \\
-u & t & 0  \\
0 & -u & t  \\
0 & 0 & -u
\end{bmatrix}$, we get  a double line with
homogeneous ideal
$$
(x_0, x_1, x_2, x_3)^2 + (x_0 t - x_1 u, x_1 t - x_2 u, x_2 t - x_3
u).
$$
\end{example}

In \cite{NNS2} the free resolution of a rope has been computed. The
key idea is to use descending induction  on the degree. The starting
point is the following observation (Lemma 3.1 in \cite{NNS2}). We
continue to use the above notation.

\begin{lemma} \label{lem-rope-ind}
Let $C$ be an $(n-k)$--rope with homogeneous ideal $I_C = ((I_L)^2, F_1, \dots, F_k)$. Then $I_{\tilde{C}} = ((I_L)^2, F_1, \dots, F_{k-1})$ defines an $(n-k+1)$--\-rope $\tilde{C}$ and the following sequence is exact
$$
0 \to I_L (-\beta_k-1) \to I_{\tilde{C}} \oplus (F_k) \to I_C \to 0.
$$
\end{lemma}

Hence, the mapping cone procedure can be used to find the minimal
free resolution of $C$, once we know the  resolution of $\tilde{C}$.
It turns out that this is most difficult when $k=1$, i.e.\ when
$I_{\tilde{C}} = I_L^2$. Then the minimal free resolution of
$\tilde{C}$ is given by an Eagon-Northcott complex. However, in
order to find the comparison maps between the resolutions of $L$ and
$\tilde{C}$, we start with the following non-minimal free
resolution. In order to state it we need more notation. We write the
Koszul complex $ P_{\bullet}$ which resolves  $ I_L $ as
\begin{equation*}\label{koszul}
P_{\bullet} \hspace*{2cm} 0 \longrightarrow \wedge^{n-1} P
\stackrel{\delta_{n-1}}{\longrightarrow} \dots \wedge^2
P \stackrel{\delta_2}{\longrightarrow} P \stackrel{\delta_1}{\longrightarrow} I_L \longrightarrow 0
\end{equation*}
where $ P = R^{n-1}(-1) = \oplus_{i=0}^r R e_i$ and $
\delta_1(e_i) = x_i $ for every i.

\begin{lemma} \label{i2res}
The ideal $(I_L)^2 $ has the following
non-minimal free resolution
$$
\begin{array}{ccccccc}
 & \wedge^{n-1} P \otimes P & & & \wedge^2 P
\otimes P \\
0  \lra & \oplus & \stackrel{\partial'_{n-1}}{\lra} & \dots & \oplus
& \stackrel{\partial'_2}{\lra} P \otimes P
\stackrel{\partial'_1}\lra (I_L)^2 \lra 0 \\[1ex]
  & \wedge^{n-1} P & & & \wedge^2 P \\
 \end{array}
$$
where
$$
\partial'_i = \left (
\begin{array}{cc}
\delta_i \otimes \id_P & (-1)^i \partial_i \\ 0 & \delta_{i}
\end{array} \right ), \quad i = 3,\ldots,n-1,
$$
$$
\partial'_2 = (\delta_2 \otimes \id_P \quad \partial_2),
\hspace*{1cm}
\partial'_1 = \delta_1 \otimes \delta_1
$$
and $ \partial_i : \wedge^i P \to \wedge^{i-1} P \otimes P $ is
the canonical map defined by $$
\partial_i (u_1 \wedge \dots \wedge u_i) = \sum_{j=1}^i (-1)^{j+1} u_1 \wedge
\dots \wedge \hat{u_j} \wedge \dots \wedge u_i \otimes u_j. $$
\end{lemma}

This is Proposition 3.2 in \cite{NNS2}. If follows by applying the mapping cone procedure to the commutative diagram
$$
\begin{array}{ccccccccc} 0 & \lra & \im
\delta_2 & \lra & I_L \otimes P & \lra & (I_L)^2 & \lra & 0 \\
 & & \uparrow  & & \uparrow \\
 & & \wedge^2 P & \stackrel{\partial_2}{\lra} & P \otimes P \\
 & & \uparrow  & & \uparrow \\
 & & \vdots & & \vdots \\
 & & \uparrow  & & \uparrow \\
 & & \wedge^{n-1} P & \stackrel{\partial_{n-1}}{\lra} & \wedge^{n-2} P
\otimes P
\\
 & & \uparrow & & \uparrow \\
 & & 0 & & \wedge^{n-1} P \otimes P \\
 & & & & \uparrow \\
 & & & & 0 \end{array} \leqno(+)
$$

We use Lemma \ref{i2res} in order to obtain a particular minimal free resolution of $I_L^2$.
Since the map $ \partial_i : \wedge^i P \to \wedge^{i-1} P \otimes P$ is split-injective we can identify $(\wedge^{i-1} P \otimes
P)/\partial_i(\wedge^i P)$ with the free $R$-module, say $D_{i-1}$, being a
direct summand of $\wedge^{i-1} P \otimes P$. Furthermore, the map $\delta_{i-1} \otimes \id_P:
\wedge^{i-1} P \otimes P \to \wedge^{i-2} P \otimes P$ induces a
homomorphism
$$
d_{i-1}: D_{i-1} \to D_{i-2}, \quad i \geq 3.
$$
Similarly, putting $D_0 = R$, the map $\delta_{1} \otimes \id_P$
induces a homomorphism $d_1: D_1 \to (I_L)^2$.

The diagram ($+$) implies that
\begin{equation*}
0 \to D_{n-1}  \stackrel{d_{n-1}}{\lra} \dots \lra D_2 \stackrel{d_2}{\lra} D_1
\stackrel{d_1}{\lra} (I_L)^2 \to 0
\end{equation*}
is an exact sequence. Comparing with the Eagon--Northcott resolution of $I_L^2$
we see that it is a minimal free resolution of $ (I_L)^2$.

Moreover, denote by
$$
\tau_B: \; Q := \oplus_{j=1}^k R(- \beta_j -1) \to P
$$
the extension of $\ffi_B$. Finally, define for $i \geq 2$
$$
\mu_{i}: \wedge^{i-1} P \otimes Q \to D_{i-1}
$$
as the composition $\wedge^{i-1} P \otimes Q \to \wedge^{i-1} P
\otimes P \to D_{i-1}$ and let $\mu_1$ be the composition $Q
\stackrel{\ffi_B}{\lra} P \stackrel{\delta_1}{\lra} R$.
\medskip

Now we can describe the minimal free resolution of an arbitrary rope
on $ L$ (cf.\ \cite{NNS2}, Theorem 3.4).

\begin{theorem} \label{Cres}
Let $ C \subset \PP^n$ be an $ (n-k)$--rope
supported on the line
$ L $ with homogeneous ideal $ I_C = ((I_L)^2, [x_0,\dots,x_r] B).$
Then we have the following
 free resolution of $ I_C$:
\begin{equation*}
G_{\bullet} \hspace*{2cm} 0 \to G_n \stackrel{d_n'}{\lra} \dots \lra G_2
\stackrel{d_2'}{\lra} G_1 \stackrel{d_1'}{\lra} I_C \to 0
\end{equation*}
where
$$
G_i := D_i \oplus (\wedge^{i-1} P \otimes Q),
$$
and
$$
d_i' := \left ( \begin{array}{cc}
d_i & (-1)^i \mu_i \\[1ex]
0 & \delta_{i-1} \otimes \id_Q
\end{array} \right )  \qquad \mif 2 \leq i \leq n-1,
$$
and
$$
D_n := 0, \qquad d_1' := \left ( \begin{array}{lr}
d_1 & -\mu_1 \
\end{array} \right ), \quad d_n' := \left ( \begin{array}{c}
(-1)^n \mu_n \\[1ex]
\delta_{n-1} \otimes \id_Q
\end{array} \right ).
$$
This resolution is minimal if and only if the rope $C \subset \PP^n$ is
non-degenerate.
\end{theorem}

\begin{remark}
Since the module $P$ has rank $n-1$, the map $d_n'$ can also be
written as $d_n' := \left (
\begin{array}{c}
(-1)^n \tau_B (-n+1) \\[1ex]
\delta_{n-1} \otimes \id_Q
\end{array} \right ).$
\end{remark}

\begin{example}[Example \ref{ex-ropes-ideal} continued] We use the
notation above.

 (i) Set $\beta := \deg F = \deg G$. Then the minimal free
 resolution of the curve $C \subset \PP^3 = \Proj (R)$ defined by $I_C = (x_0, x_1)^2 + (x_0 F + x_1
 G)$ is
 \begin{eqnarray*}
\lefteqn{
0 \lra R (-\beta - 3) \stackrel{M_2}{\lra} R (- \beta -
2)^2 \oplus R(-3) ^2 \stackrel{M_1}{\lra} } \\[1ex]
&& \hspace*{5cm}  R (-\beta - 1) \oplus R(-2)^3 \lra I_C \lra 0
\end{eqnarray*}
where
$$
M_2 = \begin{bmatrix} -F \\
-G \\
x_1 \\
-x_0
\end{bmatrix} \quad \mbox{and} \quad M_1 = \begin{bmatrix}
x_1 & 0 & F & 0 \\
-x_0 & x_1 & G & F \\
0 & -x_0 & 0 & G \\
0 & 0 & -x_0 & -x_1
\end{bmatrix}.
$$

(ii) Consider the double line $C \subset \PP^5 = \Proj (R)$ with
homogeneous ideal $I_C = (x_0, x_1, x_2, x_3)^2 + (x_0 t - x_1 u,
x_1 t - x_2 u, x_2 t - x_3 u)$. We have $P = R(-1)^4$ and $Q = R
(-2)^3$. Hence, the minimal free resolution of $C$ is of the form
\begin{eqnarray*}
\lefteqn{ 0 \lra R(-7)^3 \lra R(-6)^{12} \oplus R(-5)^4 \lra
R(-5)^{18} \oplus R(-4)^{15} \lra } \\[1ex]
& & R(-4)^{12} \oplus R(-3)^{20} \lra R(-3)^3 \oplus R(-2)^{10} \lra
I_C \lra 0.
\end{eqnarray*}
Of course, the linear part of this resolution is just the resolution
of $(x_0, x_1, x_2, x_3)^2$.
\end{example}

Let us now consider double structures of higher dimension. Manolache
has shown that often these structures are very simple.

\begin{theorem} \label{thm-lci} Let $X \subset \PP^n$ be a
  Cohen-Macaulay  double
  structure on a linear subspace of codimension $c \geq 2$. Assume
  that

\begin{tabular}{cl}
 & {\rm (i)} $c = 2$ and $n \geq 4$ \\
or & {\rm (ii)}  $X$ is locally a complete intersection and
$\frac{n+1}{3} > c$.
\end{tabular} \\
Then, after a change of coordinates, the homogeneous ideal of $X$ is
$I_X = (x_0^2, x_1,\ldots,x_{c-1})$.

\end{theorem}

\begin{proof}
These are special cases of \cite{M-92}, Theorem 1, and \\
\cite{M-95}, Theorem.
\end{proof}

The result is no longer true if the codimension of $X$ is not too small.

\begin{example}
Consider the surface $X \subset \PP^5$ with homogeneous ideal
$$
I_X := (x_0, x_1, x_2)^2 + (x_0 x_4 - x_1 x_3, x_0 x_5 - x_2 x_3, x_1
x_5 - x_2 x_4).
$$
This surface has degree two and its cohomology is
$$
H^1_* (\cI_X) \cong K \quad \mbox{and} \quad H^2_* (\cI_X) = 0.
$$
Thus, $X$ is arithmetically Buchsbaum, but not arithmetically Cohen-Macaulay.

Note, that the minimal free resolution of $X$ is linear. Hence, its
graded Betti numbers coincide with those of the disjoint union of two
planes.
\end{example}

The example shows that double structures on a linear space are not
always trivial. Nevertheless, it seems feasible to completely
characterize them, i.e.\ to extend the result about double lines to
higher extension.  Thus, we pose.

\begin{problem}
Describe the homogeneous ideals and minimal free resolutions of
Cohen-Macaulay double structures on a linear space.
\end{problem} 

For some initial results we refer to \cite{BN}, where the authors investigate the problem when a rope supported on a line can be extended to a rope that is supported on a linear space of dimension $\geq 2$. 

It is natural to consider  the next cases, i.e.\ subschemes of
degree 3, 4, \ldots. The most complicated scheme structures will
again occur for irreducible schemes that are supported on a linear
space. In order to describe these multiple structures several
filtrations have been proposed and used by B\u anic\u a and Forster
and by Manolache. We discuss them briefly.

Let $L \subset \PP^n$ be a linear subspace and let $X \subset L$ be
a Cohen-Macaulay scheme that is supported on $L$. Then we denote by
$L^{(i+1)}$ the $i$-th infinitesimal neighborhood of $L$. It is
defined by the ideal $I_L^{i+1}$.

In case $X$ is a curve, i.e.\ $L$ is a line, B\u anic\u a and
Forster introduced in \cite{BF} the Cohen-Macaulay filtration of $X$
$$
L = X_1 \subset X_2 \subset \ldots \subset X_k = X
$$
for some $k \geq 1$ where $X_i$ is obtained  by removing the
embedding points from $X \cap L^{(i)}$. The quotients
$\cI_{X_i}/\cI_{X_{i+1}}$ are vector bundles on $L$. Analyzing these
quotients, B\u anic\u a and Forster obtained an analytic
classification of the multiple structure on a line in codimension
two, up to degree $4$ (cf.\ also  \cite{BF0}). The global
classification of the
degree three structures is given  in \cite{Nollet}.  The  degree four
curves are described in 
 \cite{Nollet-Schl}. Curves of degree three and 
arbitrary codimension are studied in \cite{NNS3}. In \cite{GGP}, Theorem 5.3, it is shown that all ropes of arithmetic genus $\geq 3$ are smoothable. 

The results for codimension two curves have been used to show that the
Hilbert scheme of space curves having degree three or four are
connected (cf.\ \cite{Nollet}, \cite{Nollet-Schl}). If this true in
all higher degrees remains an open problem.

Manolache has obtained results about certain components of the Hilbert scheme of degree four space curves \cite{M-02}. Results of this type have been used to study certain moduli spaces \cite{M-81}, \cite{BM}. 

Theorem \ref{Cres}, in particular the information about the maps in
the free resolution, has been crucial for studying the component of
the Hilbert scheme that contains ropes (cf.\ \cite{NNS2}). In
particular, there is exactly one such component if the Hilbert
scheme contains ropes at all. In most cases, this component is
generically smooth.

The  B\u anic\u a-Forster  filtration in higher dimension has been
studied in \cite{Holme} and  \cite{V}.
\smallskip

Manolache has introduced two further filtrations that work in any
dimension and codimension (cf.\ \cite{M-86}, \cite{M-94}). Inspired by
liaison theory he considers the
schemes  defined by
$$
I_X : I_L^i
$$
and
$$
I_X : (I_X : I_L^i).
$$
He used his first filtration in \cite{M-92} to give a classification
of all multiple structures $X$ on $L$ when $X$ has  codimension two, its
dimension is at least two, and its degree is at most four.
In \cite{M-94} he obtained the analytic classification of multiple
structures, up to degree $6$. Combined with a result of Faltings
giving conditions on a vector bundle to guarantee that it splits as a
direct sum of line bundles, he completely described all the locally
complete intersection multiple
structures on $L$, up to degree 6, provided
the codimension of $L$ is sufficiently small. Then these structures are
all complete intersections (\cite{M-95}). For more details, we refer
to \cite{M-02}.

It would be extremely interesting to classify all Cohen-Macaulay
multiple structures on a linear space, up to degree three. A test is
the following conjecture of Vatne (\cite{V}, Conjecture 4.3), where we
slightly
deviate from our notation above.

\begin{conjecture} \label{conj-Vatne}
Let $L \subset \PP^N$ be a subspace of dimension $n \geq 6$ and
arbitrary codimension. Let $Y
\subset \PP^N$ be any Cohen-Macaulay subscheme of degree three such
that
$$
L \subset Y \subset L^{(2)}.
$$
Then there is a degree two Cohen-Macaulay subscheme $Z$ such that
$$
L \subset Z \subset Y.
$$
\end{conjecture}

This is related to the famous conjecture of Hartshorne about
codimension two varieties.

\begin{conjecture} \label{conj-Hartsh}
Every smooth variety $X \subset \PP^n$ of codimension two is a
complete intersection, provided $n \geq 6$.
\end{conjecture}

\begin{remark} \label{rem-relation}
Vatne has shown in \cite{V}, Theorem 4.4, that, for each $n$, his
conjecture is equivalent to Hartshorne's conjecture.
\end{remark}


\begin{thebibliography}{ABCD} 

\bibitem[Ballico-Notari]{BN} 
E.\ Ballico, R.\ Notari, {\em Ropes on linear subspaces of a
  projective space}, Preprint, 2005.  

\bibitem[B\u anic\u a-Forster-81]{BF0}
C.\ B\u anic\u a, O.\ Forster, {\em Sur les structures multiples},
Unpublished manuscript, 1981.  

\bibitem[B\u anic\u a-Forster-86]{BF}
C.\ B\u anic\u a, O.\ Forster, {\em Multiplicity structures on space
  curves},  Contemp.\ Math.\ {\bf 58} (1986), 47--64. 
  
\bibitem[B\u anic\u a-Manolache]{BM}
C.\ B\u anic\u a, N.\ Manolache, {\em Rank $2$ stable vector
 bundles on $\PP^3(\CC)$ with Chern classes $c_1=-1,\;c_2=4$}, 
 Math.\ Z.\ {\bf  190}  (1985), 315--339.   

\bibitem[Brodmann-Schenzel]{BS}
M.\ Brodmann, P.\ Schenzel, {\em Curves of degree $r+2$ in $\PP^r$:
  Cohomological, geometric, and homological aspects}, J.\ Algebra {\bf
  242} (2001), 577--623.

\bibitem[Bruns-Herzog]{BH}
W.\ Bruns, J.\ Herzog,
{\em Cohen-Macaulay rings.\ Rev.\ ed.},
Cambridge Studies in Advanced Mathematics {\bf 39},
Cambridge University Press,
Cambridge,
1998.

\bibitem[{Ein-Eisenbud-Katz}]{EEK}
L.\ Ein, D.\ Eisenbud, S.\ Katz, {\em Varieties cut out by quadrics:
  Scheme-theoretic versus homogeneous generation of ideals}, In:
Algebraic geometry (Sundance, UT, 1986),  51--70, Lecture Notes in
Math.\ {\bf 1311}, Springer, Berlin, 1988.

\bibitem[Eisenbud-Harris]{EH}
D.\ Eisenbud, J.\ Harris, {\em On varieties of minimal degree (A
  centennial account)}, Proceedings of Symposia in Pure Mathematics
{\bf 46} (1987), 3--13.

\bibitem[Eisenbud-Koh]{EK}
D.\ Eisenbud, J.\ Koh, {\em Some linear syzygy conjectures}, Adv.\
Math.\ {\bf 90} (1991), 47--76.

\bibitem[Fujita-1982]{F1}
T.\ Fujita, {\em Classification of projective varieties of $\Delta
  $-genus one},   Proc.\ Japan Acad.\ Ser.\ A Math.\ Sci.\ {\bf 58}
(1982), 113--116.

\bibitem[Fujita-1984]{F2}
T.\ Fujita, {\em Projective varieties of $\Delta$-genus one}, In:
Algebraic and topological theories (Kinosaki, 1984),  149--175,
Kinokuniya, Tokyo, 1986. 

\bibitem[Gallego-Gonzalez-Purnaprajna]{GGP}
F.\ J.\ Gallego, M.\ Gonzalez, B.\ P.\ Purnaprajna, {\em Smoothing of
  ropes on curves}, Preprint, 2005.  

\bibitem[Harris]{Ha}
J.\ Harris, {\em A bound on the geometric genus of propjective
  varieties}, Ann.\ Scuola Norm.\ Sup.\ Pisa Cl.\ Sci.\ (4) {\bf 8}
(1981), 36--68.

\bibitem[Hoa-St\"uckrad-Vogel]{HSV}
L.\ T.\ Hoa, J.\ St\"uckrad, W.\ Vogel, {\em Towards a structure
  theory for projective varieties of degree = codimension + 2}, J.\
Pure Appl.\ Algebra {\bf 71} (1991), 203--231.

\bibitem[Hoa]{H}
L.\ T.\ Hoa, {\em On the minimal free resolutions of projective
  varieties of degree = codimension + 2}, J.\ Pure Appl.\ Algebra {\bf
  87} (1993), 241--250.

\bibitem[Holme]{Holme}
A.\ Holme, {\em On linear subspaces with nilpotent structure},
Preprint, 1989.

\bibitem[Ionescu]{I}
P.\ Ionescu, {\em On varieties whose degree is small with respect to
codimension},  Math.\ Ann.\ {\bf 271}  (1985), 339--348. 

\bibitem[Manolache-81]{M-81}
N.\ Manolache, {\em Rank $2$ stable vector bundles on $\PP^{3}$ with
Chern classes $c_1=-1$, $c_{2}=2$},  Rev.\ Roumaine Math.\ Pures
Appl.\  {\bf 26}  (1981), 1203--1209. 

\bibitem[Manolache-86]{M-86}
N.\ Manolache, {\em Cohen-Macaulay nilpotent structures},  Rev.\
Rou\-maine Math.\ Pures Appl.\ {\bf 31}  (1986), 563--575.

\bibitem[Manolache-92]{M-92}
N.\ Manolache, {\em Codimension two linear varieties with nilpotent
  structures}, Math.\ Z.\ {\bf 210} (1992), 573--580.

\bibitem[Manolache-94]{M-94}
N.\ Manolache, {\em Multiple structures on smooth support},   Math.\
Nachr.\  {\bf 167}  (1994), 157--202.

\bibitem[Manolache-95]{M-95}
N.\ Manolache, {\em Nilpotent lci structures on global complete
  intersections}, Math.\ Z.\ {\bf 219} (1995), 403--411.

\bibitem[Manolache-03]{M-02}
N.\ Manolache, {\em Cohen-Macaulay nilpotent schemes}, In: Recent Advances in Geometry and Topology, Proceedings of the third German-Romanian Seminar on Geometry, Cluj-Napoca, 2003, 235--248. 

\bibitem[Manolache-05]{M-05}
N.\ Manolache, {\em Linkage extensions}, Preprint, 2005. 

\bibitem[Migliore]{mig}
J.\ Migliore, {\em On linking double lines}, Trans.
Amer.\ Math.\ Soc.\ {\bf 294} (1986), 177--185.

\bibitem[Nagel-Notari-Spreafico-03]{NNS1} U.\ Nagel,  R.\ Notari, M.\
  L.\ Spreafico, {\em Curves
  of degree two
and ropes on  a line: their ideals and even liaison classes}, J.\
  Algebra {\bf 265} (2003), 772-793.

\bibitem[Nagel-Notari-Spreafico-04]{NNS2} U.\ Nagel,  R.\ Notari, M.\
  L.\ Spreafico, {\em The Hilbert
scheme of degree two curves and certain ropes}, Preprint, 2004.

\bibitem[Nagel-Notari-Spreafico-05]{NNS3} U.\ Nagel,  R.\ Notari, M.\
  L.\ Spreafico, in preparation.

\bibitem[Nagel]{aBM-div}
U.\ Nagel, {\em Arithmetically Buchsbaum divisors on
    varieties of minimal degree}, Trans.\ Amer.\ Math.\ Soc.\ {\bf 351}
    (1999), 4381--4409.

\bibitem[Nollet]{Nollet}
S.\ Nollet, {\em The Hilbert scheme of degree three curves}, Ann.\
Sci.\ Ec.\ Norm.\ Sup.\ {\bf 30} (1997), 367--384.

\bibitem[Nollet-Schlesinger]{Nollet-Schl}
S.\ Nollet, E.\ Schlesinger, {\em Hilbert schemes of degree four
  curves}, Compositio Math.\ {\bf 139} (2003), 169-196. 
  
\bibitem[Okonek]{O}  
C.\ Okonek, {\em Moduli reflexiver Garben und Fl\"achen von kleinem
Grad in $\PP^{4}$},  Math.\ Z.\ {\bf 184}  (1983),  549--572. 

\bibitem[Stanley]{S} R.\ Stanley, {\em Hilbert functions of graded algebras},
Adv.\  Math.\ {\bf 28} (1978), 57--82.

\bibitem[Vatne]{V}
J.\ E.\ Vatne, {\em Multiple structures}, Preprint, 2002.


\end{thebibliography}
\end{document}